\documentclass[a4paper,12pt,reqno]{amsart}
\usepackage{amsmath,amssymb,amsthm}
\usepackage{mathrsfs}
\usepackage{enumerate}
\usepackage{graphicx}
\numberwithin{equation}{section}
\theoremstyle{plain}
\newtheorem{thm}{Theorem}[section]
\newtheorem{example}[thm]{Example}

\newtheorem{pro}{Proof}
\usepackage{mathrsfs}

\begin{document}
\title{ON CONVOLUTION FOR GENERAL NOVEL FRACTIONAL WAVELET TRANSFORM}
\author[Ashish Pathak Prabhat Yadav and M M Dixit ]{Ashish Pathak* , Prabhat Yadav** and M M Dixit** \\
 *Department of Mathematics \& Statistics\\
 Dr. Harisingh Gour Central University \\
    Sagar-470003, India.\\.
    ** Department of Mathematics \\
NERIST, Nirjuli-791-109, India\\}
\date{}
\thanks{$^{*}$E-mail: pathak\_maths@yahoo.com}
\begin{abstract}
Using Pathak and  Pathak techniques,  the basic function $D^\alpha(u,v,w)$ associated with general novel fractional wavelet transform (GNFrWT)is defined and its properties are investigated. By using basic function $D^\alpha(u,v,w)$ translation and convolution associated with GNFrWT are defined and certain existence theorems are proved for basic function  and associated convolution.
\end{abstract}
\maketitle
\section{Introduction} The general novel fractional wavelet transfrom (GNFrWT) of a function $h(t)$ with respect to wavelet $\phi$ can be defined as [1],
\begin{equation}
\label{1.1}
 (W^\alpha_\phi h)(a,b) = \int_{\mathbb{R}} h(t)
 \,\,    \overline{\phi^{\alpha}_{a,b}(t)}\,dt,
\end{equation}
where
\begin{eqnarray}
\label{1.2}
\phi^{\alpha}_{a,b}(t) = e^{(-i/2)(t^2-b^2)cot\theta} |a|^{-\rho} \phi{\left(\frac{t-b}{a}\right)};\,\, a\in{\mathbb{R_+}}\,,\, b\in{\mathbb{R}}\,and \,\rho\geq 0.
 \end {eqnarray}
The inversion formula  for (\ref{1.1}) with respect to the $e^{(-i/2)(t^2-b^2)cot\theta} \,\phi_{a,b}(t)$ can be given as
\begin{eqnarray}
\label{1.3}
h(t)= \frac{1}{c_\phi} \int_{\mathbb{R}}\int_{\mathbb{R}}(W^\alpha_\phi h)(a,b) \phi^{\alpha} _{a,b}(t) |a|^{2\rho -3}dbda ,
\end{eqnarray}
where
\begin{eqnarray}
\label{1.4}
O<C_\phi    =  \displaystyle\int_{\mathbb{R}}\frac{|\hat\phi(\omega)|^2}{|\omega|}\,\,d\omega<\infty.
\end{eqnarray}
 If $\alpha=1 $, the GNFrWT , correlate with the wavelet transform (WT).
 As \cite{shi} , the  GNFrWT can be expressed as in terms of the fractional fourier transform $H^\alpha(v)$ of the signal $h(t)$.
 \begin{equation}
 \label{1.5}
 (W^\alpha_h)(a,b) = \int_{\mathbb{R}}\sqrt{2\pi}\,\, a^{-\rho+1} H^\alpha(v)\overline{\hat\phi(avcosec\theta)}\,\,
\overline {K_ {-\alpha}(v,b)}dv
\end{equation}
where $ \hat\phi(a v cosec\theta)$ indicates the fourier transform of  $ \phi(t)$.
Also, the GNFrWT can be rewritten as \cite{shi}
\begin{eqnarray}
\label{1.6}
(W^\alpha_\phi h)(a,b) = e^{(-i/2){b^2}cot\theta}\int_{\mathbb{R}}h(t)e^{(i/2){t^2}cot\theta}\overline{\phi_{a,b}(t)}dt.
\end{eqnarray}
\begin{thm}{( Parseval formula).}
 If wavelet $ \phi\in L^2(\mathbb{R}) $  and $ (W^\alpha_\phi h)(a,b)$ is the GNFrWT of $ h\in L^2(\mathbb{R}) $ , then for any function $ h,g \in L^2({\mathbb{R}})$ ,
 \begin{eqnarray}
 \label{1.7}
 \int _{\mathbb{R}}\int _{\mathbb{R}}(W^\alpha_\phi h){(a,b)}\,\, \overline{(W^\alpha_\phi g){(a,b)}}\,\, |a|^{2\rho-3}dbda = C_\phi(h,g).
 \end{eqnarray}                                                                                                           where
 \begin{eqnarray}
 \label{1.8}
  O<C_\phi  =  \nonumber \displaystyle\int_{\mathbb{R}}\frac{|\hat\phi(\omega)|^2}{|\omega|}\,\,d\omega<\infty.
  \end{eqnarray}
\end{thm}
 \begin{pro}
 By using (\ref{1.5}),we have
 \begin{eqnarray}
 \label{1.9}
 W^\alpha_ \phi (h)(a,b) = \int_{\mathbb{R}}\sqrt{2\pi} \, \, {a^{{-\rho + 1}}} {H^\alpha (h)(u)}\overline{ \hat\phi(aucosec\theta)}\,\,\overline{K_ {-\alpha}(u,b)} du
 \end{eqnarray}
 \begin{eqnarray}
 \label{1.10}
 \overline{W^\alpha_ \phi (g)(a,b)}= \int_{\mathbb{R}} \sqrt {2\pi} \,\, { a^{{-\rho + 1}}}\overline{H^\alpha (g)(v)}\,\,{  \hat \phi(avcosec\theta)} {K_{-\alpha}(v,b)} dv
 \end{eqnarray}
 from (\ref{1.9}) and (\ref{1.10}), we get
  \begin{eqnarray}
  \label{1.11}
  \int_{\mathbb{R}}\int_{\mathbb{R}} W^\alpha_ \phi (h)(a,b)\, \overline{W^\alpha_ \phi (g)(a,b) }\,\ |a|^{2\rho-3} dbda = <h,g> C_\phi.
\end{eqnarray}
\end{pro}
 \begin{thm} (Inversion Formula ) If wavelet $ \phi\in L^2(\mathbb{R}) $  and $ (W^\alpha_\phi h)(a,b)$ is the GNFrWT of $ h\in L^2(\mathbb{R}) $ , then the reconstruction of \, $ h $ is given by
 \begin{eqnarray}
  \label{1.12}
   h(u) =    \frac{1}{C_\phi} \int_{\mathbb{R}} \int_{\mathbb{R}} (W^\alpha _\phi h)(a,b) \phi^{\alpha} _{a,b} {(t)}  \,\,|a|^{2\rho-3} dbda
\end{eqnarray}
\end{thm}
\begin{pro}
 By using (\ref{1.11}) for $ h = g $ , then we can find (\ref{1.12}).\\
 This paper is arranged in the  following manner :- In the next section, we define basic function $ D^\alpha(u,v,w) $, translation and associated convolution  for GNFrWT.  In the last third section, we obtained and established the certain existence theorem and convolution theorem , by using Pathak and Pathak techniques \cite{Pathak}

\section{ Basic function , translation  and  associated convolution  for GNFrWT }
Now, by  using  Pathak and Pathak  techniques \cite{Pathak}, we define the basic function $D^\alpha(u,v,w)$ , translation ${\tau^\alpha_u}$ and associated  convolution $\#^\alpha $ operators for GNFrWT.\\
 The basic function \,$ D^\alpha(u,v,w)$ for (\ref{1.1}) is define as
\begin{eqnarray}
\label{2.1}
W^\alpha_\phi[D^\alpha(u,v,w)](a,b) &=&\int_{\mathbb{R}} D^\alpha(u,v,w)\overline{{\phi^\alpha}_{a,b}(t)}dt\nonumber \\
 &=& \overline{\psi^\alpha_{a,b}(w)}\,\ \overline{\chi^\alpha_{a,b}(v)},
\end{eqnarray}
where $ \psi^\alpha , \phi^\alpha $ and $ \chi^\alpha $ are three  fractional wavelets satisfying certain conditions (\ref{1.2}).\\
Now, by  using(\ref{1.3}) we get,
\begin{equation}
\label{2.2}
D^\alpha(u,v,w)= C^{-1}_\phi \int_{\mathbb{R}}\int_{\mathbb{R}}\overline{\psi^\alpha_{a,b}(w)}\,\ \overline{\chi^\alpha_{a,b}(v)}\ {\phi^\alpha_{a,b}(u)} \, |a|^{2\rho-3}da db.
\end{equation}
The translation  $\tau^\alpha_u $ is defined as \cite{Pathak}
\begin{eqnarray*}
{(\tau^\alpha_u h)(v)} & = & h^*(u,v)=\int_{\mathbb{R}} D^\alpha(u,v,w)h(w)dw
\\ &=& C^{-1}_\phi\int_{\mathbb{R}}\int_{\mathbb{R}}\int_{\mathbb{R}}\overline{\psi^\alpha_{a,b} (w)}\,\, \overline{\chi^\alpha_{a,b}(v)}\,\, \phi^\alpha_{a,b}(u)\, h(w)|a|^{2\rho-3}dadbdw.
\end{eqnarray*}
The associated  convolution  is defined  as
\begin{eqnarray}
\label{2.3}
(h\#^\alpha g)(u)\nonumber &=& \int_{\mathbb{R}}h^*(u,v) g(v)dv \\ \nonumber&=& \int_{\mathbb{R}}\int_{\mathbb{R}}D^\alpha(u,v,w)\,h(w)\,g(v)dv dw\\ \nonumber &=& C^{-1}_\phi \int_{\mathbb{R}}\int_{\mathbb{R}}\int_{\mathbb{R}}\int_{\mathbb{R}} \overline{{\psi^\alpha}_{a,b}(w)}\,\,  \overline{{\chi^\alpha}_{a,b}(v)} \,\,\phi^\alpha_{a,b}(u) h(w)g(v)\left|a\right|^{2\rho-3}dadbdwdv. \\
\end{eqnarray}
\end{pro}
\begin{example}
\textbf{Basic function $ D^\alpha(u,v,w) $  for general novel fractional morlet wavelet transform}
\\Let $\psi(t) = \chi(t)=\phi(t) =   e^{iw_o t - \frac{1}{2}{t^{2}}} $ be a morlet wavelet \cite{debnath} . Then  the general novel fractional morlet wavelet transform is given by $
\phi^\alpha _{a,b}(t) = |a|^{-\rho} e^{{iw_o {\frac{(t-b)}{a}}} - {\frac{1}{2}{{[{\frac{(t-b)}{a}}]}}^{2}}}$ .
Now, from (\ref{2.2})
\begin{eqnarray*}
D^\alpha(u,v,w)
& = &
{C^{-1}_\phi}\int_{\mathbb{R}}\int_{\mathbb{R}} \,e^{{(i/2)}({w^2+v^2-u^2-b^2)cot\theta}} {\overline{\psi_{a,b}(w)}}\,\,{\overline{\chi_{a,b}(v)}}\,\,\phi_{a,b}(u)\left|a\right|^{2\rho-3}dadb.
\\ &=&
 {C^{-1}_{\phi}}e^{{(i/2)}{(w^2+v^2-u^2)}cot\theta}\int_{\mathbb{R}}\int_{\mathbb{R}}
e^{{(-i/2)}{b^2}cot\theta} \,\,
e^{{-iw_o}{{\frac{(w-b)}{a}}}-{\frac{1}{2}}{\frac{(w-b)^2}{a^2}}}
 \\ &  & \times
e^{{-iw_o}{{\frac{(v-b)}{a} }} -  {\frac{1}{2}} {\frac{(v-b)^2}{a^2}}}\,\,
e^{{ iw_o}{\frac{(u-b)}{a}}-{\frac{1}{2}} {\frac{(u-b)^2}{a^2}}} \,\,|a|^{\rho-3} dbda
\\ &=&
 {C^{-1}_{\phi}}e^{{(i/2)} {(w^2+v^2-u^2)}cot\theta}
                       \int_{\mathbb{R}}\int_{\mathbb{R}}\,\, e^{{(-i/2)}{b^2}cot\theta}
\\ &  & \times
                        e^{(itw_o)(b+u-w-v)} \,\,e^{-t^2 \frac{(w-b)^2+(v-b)^2+(u-b)^2}{2}} \,\,|t|^{\rho+1}dtdb
  \\ &=&
 2 {C^{-1}_{\phi}}e^{{(i/2)}{(w^2+v^2-u^2)}cot\theta}
                            \int_{\mathbb{R}}\int_{\mathbb{R}} \,\,e^{{(-i/2)}{b^2}cot\theta}  \cos{[w_o t(b+u-w-v)]}
\\ & & \times
                            e^{-t^2 \frac{(w-b)^2+(v-b)^2+(u-b)^2}{2}} \,\,|t|^{\rho+1} dtdb
\\ & =&
{C^{-1}_{\phi}}e^{{(i/2)}{(w^2+v^2-u^2)}cot\theta}\,\,\Gamma{(1+\rho/2)
                        {2^{(1+\rho/2)}}}
\\ & & \times
                        \int_{\mathbb{R}}
                          [{(w-b)^2 +(v-b)^2 +(u-b)^2}]^ {-1-\rho/2}
\\ &  & \times
 e^{{(-i/2)}{b^2}cot\theta} {_1}F_{1} {(1+\rho/2;1/2};\,\,
                         \frac{- w^2_o (b+u-w-v)^2}{ 2[(w-b)^2 +(v-b)^2 +(u-b)^2]}db , \rho>0 ,
 \end{eqnarray*}
by [\, \cite{erdelyi},p.15(14)] where $ _1F_{1}( a;b;u)$ is confluent hypergeometric function.
\end{example}
\begin{example}\textbf{Basic function $D^\alpha(u,v,w)$  for  general novel fractional\\ mexican hat wavelet tansform } \\
 The corresponding  Mexican-Hat wavelet \cite{debnath} is $ \psi(t) = \chi(t)=\phi(t)=  (1-t^2)\,\, e^{\frac{- t ^2}{2}}$ . Then  the general novel fractional mexican hat wavelet transform is given by $\phi^\alpha_{b,a}{(t)}={|a|^{-\rho}}{\left(1-\frac{(t-b) ^2}{a ^2}\right)} \,\,e^ {{\frac{-(t-b)^2}{2 {a}^{2}}}}$ .\\
Now by using (\ref{2.2}), we have
\begin{eqnarray*}
D^\alpha(u,v,w)
   & = & C^{-1}_{\phi} {e^{(i/2)(w^2+v^2-u^2)cot\theta}}\int_{\mathbb{R}} \int_{\mathbb{R}} {e^{(-i/2)b^2 cot\theta}}
  \left ( 1-\frac{(w-b)^{2}}{a^2} \right) \\ && \times  \left( 1-\frac{(v-b)^2}{a^2}\right)  \left( 1-\frac{(u-b)^2}{a^2}\right)\\ && \times e^{(- 1/2a^2) ((w-b)^2 +(v-b)^2 +(u-b)^2)}  dbda.
 \\    & = &
    C^{-1}_{\phi} e^{(i/2)(w^2+v^2-u^2)cot\theta} \int_{\mathbb{R}} e^{(-i/2)b^2 cot\theta} e^{(-t^{2} {L/2})}\,\,t^{\rho+1}\\  && \times (t^{4} - N t^{2} +M) dt db.
  \\ & = &
  {C^{-1}_{\phi}} {e^{(i/2)(w^2+v^2-u^2)cot\theta}}\int_{\mathbb{R}}{e^{(-i/2)b^2 cot\theta}}db  {(-1/2)(-1)^{\rho + 1}}
  \\ & & \times
    \bigg( \Gamma((\rho+6)/2) (L/2)^{-(\rho+6)/2} -\Gamma((\rho+4)/2) N(L/2)^{-(\rho+4)/2}  \\ & & + \Gamma((\rho+2)/2) M (L/2)^{\frac{-(\rho+2)}{2}} \bigg)
   \end{eqnarray*}
by   [\, \cite{erdelyi},p.313(13)] , where $ \rho> 0, \,\,  L=[(w-b)^2 +(v-b)^2+(u-b)^2],\\ M= [{(u-b)^2 (v-b)^2} + {(u-b)^2(w-b)^2 }+{(v-b)^2(w-b)^2}] and\\
N = [L +(u-b)^2(v-b)^2(w-b)^2]$.
\end{example}
In the following section  we have obtained boundedness result for the basic function $D^\alpha(u,v,w)$ and then establish existence theorem for the general novel fractional wavelet convolution  and prove  $ W^\alpha_\phi(h\#^\alpha g)= (W^\alpha_\psi{h})(W^\alpha_\chi{g}).$
\section{ Existence Theorems}
First we obtain boundedness results for the basic function $D^\alpha(u,v,w)$.
\begin{thm}
\label{3.1}
Let $ (1+\left|z\right|^\rho)\phi(z)\in L^{p}(\mathbb{R}),\chi\in L^q(\mathbb{R}), \frac{1}{p}+\frac{1}{q}=1 $ and \\ $(1+\left|z\right|^\rho)\psi(z)\in L^{1}(\mathbb{R}),\rho\geq 0 $.Then
\begin{eqnarray}
\label{3.2}
\left|D^\alpha(u,v,w)\right| & \leq\nonumber & 2^{\rho+\frac{1}{p}}C^{-1}_\phi \left|v-w\right|^{-\frac{1}{q}}\left|u-w\right|^{-\frac{1}{p} -\rho }\left\|\chi\right\|_q\left\|(1+\left|z\right|^\rho)\phi(z)\right\|_p \\&&
\times \; \left\|(1+\left|z\right|^\rho)\psi(z)\right\|_1,
\end{eqnarray}
 where $ O<C_\phi  =\displaystyle\int_{\mathbb{R}}\frac{|\hat\phi(\omega)|^2}{|\omega|}\,\,d\omega<\infty.$\\
\end{thm}
\begin{pro}
From (\ref{2.2}), we have
\begin{eqnarray}
\label{3.3}
|D^\alpha(u,v,w)| &=&  \nonumber   |C^{-1}_\phi \int_{\mathbb{R}}\int_{\mathbb{R}}
                            |a|^{-\rho}\,\overline{e^{{(-i/2)(w^2-b^2)cot\theta}}}
                            \,\, \overline{ \psi \left  (\frac{(w-b)}{a} \right )}\,\,|a|^{-\rho}
                            \overline{e^{{(-i/2)(v^2-b^2)cot\theta}}}
                                         \\ &  &  \times \nonumber \overline{ \chi\left  ( \frac{(v-b)}{a}\right)}
                             |a|^{-\rho} \,\ e^{{(-i/2)(u^2-b^2)cot\theta}} \,\,
                                  \,\,{ \phi \left ( \frac{(u-b)}{a}\right)} \,\,|a|^{2\rho-3} db da |
\\ &=& \nonumber
C^{-1}_\phi \int_{\mathbb{R}}\int_{\mathbb{R}}  |\overline{\psi \left (\frac{(w-b)}{a}\right)} | \, | \, \overline{\chi \left( \frac{(v-b)}{a}\right ) }\, |
      \,\,  |\,{ \phi\left ( \frac{(u-b)}{a} \right)}  |\,\,|a|^{- \rho-3} db da \\
\end{eqnarray}
 By using Theorem 2.1 \cite{Pathak} we get the required result.
\end{pro}
\begin{thm}
\label{3.4}
\,\, (i) \,\,Let $ \psi \in L^1(\mathbb{R}), \phi \in L^p(\mathbb{R}), \chi \in L^q(\mathbb{R}),
p,q,>1, 0<\rho<1 $ and $\frac{1}{p}+ \frac{1}{q}=1+\rho$.Then
\begin{eqnarray}
\label{3.5}
\int_\mathbb{R}\left|D^\alpha(u,v,w)\right|dw \leq C^{-1}_\phi   C(p,\rho)\left|u-v\right|^{-\rho}\left\|\psi\right\|_1\left\|\phi\right\|_p\left\|\chi\right\|_q,
\end{eqnarray}
where $C( p,\rho)$ is a constant.\\
(ii)\,\ Let $ \psi \in L^1(\mathbb{R}), (1+\left|y\right|^{\rho-1})\phi\mathfrak{}(y) \in L^1(\mathbb{R}), (1+\left|y\right|^{\rho-1})\chi(v) \in L^1(\mathbb{R})$,and \\$ \rho \geq 1$.Then
\begin{eqnarray}
\label{3.6}
\int_\mathbb{R}\left|D^\alpha(u,v,w)\right|dw \nonumber & \leq &  C^{-1}_\phi 2^{\rho-1}\left|u-v\right|^{-\rho}\left[\left\|\phi(x)x^{\rho-1}\right\|_1\left\|\chi\right\|_1+  \left\|\chi(y)y^{\rho-1}\right\|_1\left\|\phi\right\|_1\right]\\ && \times \; \left\|\psi\right\|_1.
\end{eqnarray}
\end{thm}
\begin{pro} From (\ref{3.3})  and using Theorem 2.2  \cite{Pathak} the required results follows.
\end {pro}
\begin{thm}
\label{3.7}
(i).Let $ \psi \in L^p(\mathbb{R}), \phi \in L^1(\mathbb{R}), \chi \in L^q(\mathbb{R}),
p,q,>1, 0<\rho<1 $ and $\frac{1}{p}+ \frac{1}{q}=1+\rho$.Then
\begin{eqnarray}
\label{3.8}
\int_\mathbb{R}\left|D^\alpha(u,v,w)\right|du \leq C^{-1}_\phi C(p,\rho)\left|v-w\right|^{-\rho}\left\|\psi\right\|_p\left\|\phi\right\|_1\left\|\chi\right\|_q.
\end{eqnarray}
(ii).Let $ \phi \in L^1(\mathbb{R}), (1+\left|x\right|^{\rho-1})\chi(x) \in L^1(\mathbb{R}), (1+\left|x\right|^{\rho-1})\psi(x) \in L^1(\mathbb{R})$,and \\$ \rho \geq 1$.Then
\begin{eqnarray}
\label{3.9}
\int_\mathbb{R}\left|D^\alpha(u,v,w)\right|du \nonumber & \leq &  C^{-1}_\phi 2^{\rho-1}\left|v-w\right|^{-\rho}\left[\left\|\psi(x)x^{\rho-1}\right\|_1\left\|\chi\right\|_1+  \left\|\chi(x)x^{\rho-1}\right\|_1\left\|\psi\right\|_1\right]\\ && \times \; \left\|\phi\right\|_1.
\end{eqnarray}
\end{thm}
The proof is similar to that Theorem 3.2 .
\begin{thm}
\label{3.10}
(i).Let $ \psi \in L^q(\mathbb{R}), \phi \in L^p(\mathbb{R}), \chi \in L^1(\mathbb{R}),
p,q,>1, 0<\rho<1 $ and $\frac{1}{p}+ \frac{1}{q}=1+\rho$.Then
\begin{eqnarray}
\label{3.11}
\int_\mathbb{R}\left|D^\alpha(u,v,w)\right|dv \leq C^{-1}_\phi C(p,\rho)\left|u-w\right|^{-\rho}\left\|\psi\right\|_q\left\|\phi\right\|_p\left\|\chi\right\|_1.
\end{eqnarray}
(ii).Let $ \chi \in L^1(\mathbb{R}), (1+\left|x\right|^{\rho-1})\phi(x) \in L^1(\mathbb{R}), (1+\left|x\right|^{\rho-1})\psi(x) \in L^1(\mathbb{R})$,and \\ $ \rho \geq 1$.Then
\begin{eqnarray}
\label{3.12}
\int_\mathbb{R}\left|D^\alpha(u,v,w)\right|dv \nonumber & \leq &  C^{-1}_\phi
2^{\rho-1}\left|u-w\right|^{-\rho}\left[\left\|\psi(x)x^{\rho-1}\right\|_1\left\|\phi\right\|_1+  \left\|\phi(x)x^{\rho-1}\right\|_1\left\|\psi\right\|_1\right]\\ && \times \; \left\|\chi\right\|_1.
\end{eqnarray}
\end{thm}
The proof is similar to that Theorem 3.2
\begin{thm}
\label{3.13}
Let $\phi \in L^1(\mathbb{R}),\psi \in L^p(\mathbb{R}),\chi \in L^q(\mathbb{R}), p,q>1, 0<\rho<1,\\\frac{1}{p}+\frac{1}{q}=\rho+1,  h \in L^{r}(\mathbb{R})$ and $ g \in L^{r'}(\mathbb{R}), r,r'>1,\frac{1}{r}+\frac{1}{r'}+\rho=2 $.Then \\
\begin{eqnarray*}
||( h\#^\alpha g)||_1 \leq C^{-1}_\phi C(\rho,p,r)\,|||\phi|_1 ||\psi||_p||\chi||_q||g||_{r'} ||h||_{r}
\end{eqnarray*}
where  $C_\phi $ is given by  (1.4)  and $ C\,(\rho,p,r)$ is a constant.
\end{thm}
\begin{pro}  we have
 \begin{eqnarray*}
  \int_\mathbb{R}{|(h\#^\alpha g)(u)|}du & \leq &\int_\mathbb{R}\left(\int_\mathbb{R}|h^{*}(u,v)| \,|g(v)|dv \right)du
\\ &=&  \int_\mathbb{R}{|g(v)|dv} \int_\mathbb{R}{|h^{*}(u,v)|}du
\\ & \leq &  \int_\mathbb{R}{|g(v)|dv} \int_\mathbb{R}\left( \int_\mathbb{R}|D^\alpha(u,v,w)|\,|h(w)|dw \right)du
\end{eqnarray*}
by using above Theorem 3.3, we get
\begin{eqnarray*}
\int_\mathbb{R}{|(h\#^\alpha g)(u)|du} \leq C^{-1}_\phi\,C(p,\rho)\,||\phi||_1\,||\chi||_q\,||\psi||_p   \int_\mathbb{R}{|g(v)|dv}\,\int_\mathbb{R}{|h(w)|}\,{|v-w|^{-\rho}}\, dw
\end{eqnarray*}
Therefore, by using Theorem 2.5 \cite{Pathak} , we get the required result.
\end{pro}
\begin{thm}
\label{3.14}
Let  $ \phi\in L^1(\mathbb{R}) \cap L^\infty(\mathbb{R}), \chi \in L^q(\mathbb{R}),\psi \in L^{p}(\mathbb{R}), p,q > 1, 0 < \rho < 1 ,  \frac{1}{p}+\frac{1}{q}   =  \rho+1 $ .Assume further that $ h \in L^r(\mathbb{R}),g \in L^{r'}(\mathbb{R}), r,r' > 1 $   and  $\frac{1}{r}+\frac{1}{r'}+\rho  =  2$ . Then
\begin{eqnarray*}
W^\alpha_\phi(h\#^\alpha g)(a,b) =  (W^\alpha_\psi h) (a,b) (W^\alpha_\chi g)(a,b)
\end{eqnarray*}
\end{thm}
\begin{pro} By using  above Theorem 3.5 , $ (h\#^\alpha g)\in\,L^{1}{(R)}$ .As  given basic wavelet $\phi\in L^{1}{(R)}\cap L^\infty{(R)},  W^\alpha_\phi(h\#^\alpha g)(a,b)$ exist.\\
 From (\ref{2.1}), we have
 \begin{eqnarray*}
 W^\alpha_\phi(h\#^\alpha g)(a,b)   & = &   \int_\mathbb{R}{(h\#^\alpha g)(u)} \overline{\phi^\alpha \left(\frac{u-b}{a}\right)} \, |a|^{-\rho} du
 \\ & = & \int_{\mathbb{R}} \overline{\phi^\alpha\left(\frac{u-b}{a}\right)} du \int_\mathbb{R}\int_\mathbb{R} D^\alpha{(u,v,w)}{h(w)}{g(v)}\,\, |a|^{-\rho} dw dv
  \\ & = &
  \int_\mathbb{R}\int_\mathbb{R} h{(w)}g{(v)}dw dv \int_\mathbb{R} D^\alpha{(u,v,w)} \overline{\phi^\alpha\left(\frac{u-b}{a}\right)}\,\, |a|^{-\rho} du
  \\ & = &
  \int_\mathbb{R}\int_\mathbb{R} h{(w)}g{(v)}dw dv \,\, \overline{\psi^\alpha_{a,b}{(w)}}\,\, \overline{\chi^\alpha_{a,b}{(v)}}
  \\ & =&
   \int_\mathbb{R} h{(w)} \overline{\psi^\alpha_{a,b}{(w)}}  dw \int_\mathbb{R} g{(v)}\overline{\chi^\alpha_{a,b}{(v)}} dv
  \\ & =&
   (W^\alpha_\psi h)(a,b)  (W^\alpha_\chi g)(a,b)
 \end{eqnarray*}
\end{pro}
\section*{Acknowledgment}
The work of the first author was supported by U.G.C start-up grant for newly recruited faculty.


\thebibliography{00}

 \bibitem {shi}  J. Shi , N.T.Zhang and X.P.Liu.  A novel fractional wavelets transform and its application; science china information science , (2012), 55(6):  1270-1279.
\bibitem{debnath} L. Debnath .  Wavelet transforms and their application ; Birkhauser, (2002) .
\bibitem{bhatnagar} G. Bhatnagar ;  Q. M. Jonathanwu and B. S. Raman .   Discrete fractional wavelet transform and its application to multiple encryption ;  information science ,  (2013), 223: 297-316.
\bibitem{erdelyi} A. Erde'lyi,  W. Magnus, F. Oberhettinger and F. G. Tricomi.  Tables of Integral Transforms, Vol. 1. McGraw-Hill, New York (1954).
\bibitem{Perrier}  V. Perrier and C. Basdevant.  Besov norms in terms of the continuous wavelet transform, application to structure function, math, models and methods in appl. Sci. ,(1996), 6(5): 649-664.
\bibitem{Pathak} R. S. Pathak and Ashish Pathak . On convolution for wavelets transform ;international journal of wavelets, multiresolution and information processing , (2008), 6(5): 739-747.
\bibitem{R.S.Pathak} R. S. Pathak.  Convolution for the discrete wavelet transform ;international journal of wavelets, multi resolution and information processing, (2011),9(6): 905-922.
\bibitem{shantha} R. Shantha,  Selva Kumari,  R .Suriya Prabha  and V.Sadasivam. ECG signal coding using biorthogonal wavelet- based burrows-wheeler coder; international journal of wavelets, multiresolution and information processing ,(2011),  9(2) : 269-281.
\bibitem{michel} Michel alves Lacerda , Rodrigo Capobianco Guido , Leonardo Mendes De Souza , Paulo Ricardo Franchi Zulato , Jussara Ribeiro , and  Shi-huang Chen.  A wavelet-based speaker verification algorithm ;international journal of wavelets, multiresolution and information processing, (2010),8(6): 905-912.
\end{document}